\documentclass[12pt]{article}
\newtheorem{definition}{Definition}
\newtheorem{proposition}{Proposition}
\newtheorem{theorem}{Theorem}

\newcommand{\field}[1]{\mathbb{#1}}
\usepackage{graphics,amsmath,amssymb}
\usepackage{amssymb}
\usepackage{natbib}
\begin{document}
\title{\Large{Are Deterministic Descriptions And Indeterministic Descriptions Observationally Equivalent?}}
\author{Charlotte Werndl\\
\normalsize The Queen's College, Oxford University, charlotte.werndl@queens.ox.ac.uk}
\date{\normalsize This is a pre-copyedited, author-produced PDF of an article accepted for publication in Studies in History and Philosophy of Modern Physics following peer review. The definitive publisher-authenticated version ``C. Werndl (2009), Are Deterministic Descriptions and Indeterministic Descriptions Observationally Equivalent, Studies in History and Philosophy of Modern Physics 40, 232-242'' is available online at: http://www.sciencedirect.com/science/article/pii/S135521980900032X}
\maketitle
\begin{abstract}
The central question of this paper is: are deterministic and indeterministic descriptions observationally equivalent in the sense that they give the same predictions? I tackle this question for measure-theoretic deterministic systems and stochastic processes, both of which are ubiquitous in science. I first show that for many measure-theoretic deterministic systems there is a stochastic process which is observationally equivalent to the deterministic system. Conversely, I show that for all stochastic processes there is a measure-theoretic deterministic system which is observationally equivalent to the stochastic process. Still, one might guess that the measure-theoretic deterministic systems which are observationally equivalent to stochastic processes used in science do not include any deterministic systems used in science. I argue that this is not so because deterministic systems used in science even give rise to Bernoulli processes. Despite this, one might guess that measure-theoretic deterministic systems used in science cannot give the same predictions at every observation level as stochastic processes used in science. By proving results in ergodic theory, I show that also this guess is misguided: there are several deterministic systems used in science which give the same predictions at every observation level as Markov processes. All these results show that measure-theoretic deterministic systems and stochastic processes are observationally equivalent more often than one might perhaps expect. Furthermore, I criticise the claims of the previous philosophy papers Suppes (1993, 1999), Suppes and de Barros (1996) and Winnie (1998) on observational equivalence.
\end{abstract}
\tableofcontents
\newpage
\section{Introduction}
There has been a lot of philosophical debate about the question of
whether the world is deterministic or indeterministic. Within this
context, there is often the implicit belief (cf.\
\citeauthor{Weingartner1996}, 1996, p.~203) that deterministic and
indeterministic descriptions are not observationally equivalent. However, the question whether these descriptions are observationally equivalent has
hardly been discussed.

This paper aims to contribute to fill this gap. Namely, the central questions of this paper are the following:
\textit{are deterministic mathematical descriptions
and indeterministic mathematical descriptions observationally
equivalent? And what is the philosophical significance of the
various results on observational equivalence?}

The deterministic and indeterministic descriptions of concern in
this paper are measure-theoretic deterministic systems and
stochastic processes, respectively. Both are ubiquitous in
science. Because of lack of space, I concentrate on descriptions
where the time varies in discrete steps; but I point out that
analogous results also hold for a continuous time
parameter.

More specifically, when saying that a deterministic system and a
stochastic process are \textit{observationally equivalent}, I mean
the following: the deterministic system, when observed, gives the
same predictions as the stochastic process. In what follows, when
I say that a stochastic process can be \textit{replaced} by a
deterministic system, or conversely, I mean that it can be
replaced by such a system in the sense that they are observationally
equivalent.

This paper proceeds as follows. In section \ref{gesund} I will introduce stochastic processes and measure-theoretic deterministic systems. In section \ref{BI} I
will show that measure-theoretic deterministic systems and stochastic processes can often be replaced by each other. Given this, one
might still guess that it is impossible to replace stochastic processes of the kinds in fact used in science by measure-theoretic deterministic systems that are used in science. One might also guess that it is impossible to replace measure-theoretic deterministic systems of the kinds used in science at every observation level by stochastic processes that are used in science.
By proving results in ergodic theory, I will show in section \ref{AI} that these two guesses are wrong.
Therefore, even kinds of stochastic processes and deterministic systems which seem to give very different predictions are observationally equivalent. Furthermore, I will criticise the claims of the previous philosophical papers Suppes (1993,1999), \cite{Suppes1996} and \cite{Winnie1998} on observational equivalence.

For a less technical treatment of the issues discussed in this paper, see \citeauthor{Werndl2009c} (\citeyear{Werndl2009c}).

\section{Stochastic processes and deterministic systems}\label{gesund}

The indeterministic and deterministic descriptions I deal with are stochastic processes and measure-theoretic deterministic systems, respectively.
There are two types of them: either the
time parameter is discrete (discrete processes and systems) or
there is a continuous time parameter (continuous processes and
systems). I consider only discrete descriptions, but analogous
results hold for continuous descriptions, and these results are
discussed in \cite{Werndl2009d}.

\subsection{Stochastic processes}\label{SP}
A stochastic process is a process governed by probabilistic laws.
Hence there is usually indeterminism in the time-evolution: if the
process yields a specific outcome, there are different outcomes
that might follow; and a probability distribution measures the
likelihood of them.
I call a sequence which describes a possible
time-evolution of the stochastic process a realisation. Nearly
all, but not all, the indeterministic descriptions in science
are stochastic processes.\footnote{For instance, Norton's dome (which satisfies Newton's laws) is indeterministic because the time evolution fails to be bijective. Nothing in Newtonian mechanics requires us to assign a probability measure on
the possible states of this system. It is possible to assign a
probability measure, but the question is whether it is natural
\cite[cf.][pp.~8--9]{Norton2003}.}

Let me formally define stochastic processes.\footnote{I assume basic knowledge about
measure theory and modern probability theory. For more details, see \cite{Doob1953}, \cite{Cornfeldetal1982} and
\cite{Petersen1983}.} A \textit{random variable} is a measurable function
$Z:\Omega\rightarrow\bar{M}$ from a probability space, i.e.~a
measure space $(\Omega,\Sigma_{\Omega},\nu)$ with $\nu(\Omega)=1$,
to a measurable space $(\bar{M},\Sigma_{\bar{M}})$ where
$\Sigma_{\bar{M}}$ denotes a $\sigma$-algebra on
$\bar{M}$.\footnote{For simplicity, I assume that any measure
space is complete, i.e.~every subset of a measurable set of
measure zero is measurable.} The probability measure
$P_{Z}(A)=P\{Z\in A\}:=\nu(Z^{-1}(A))$ for all $A\in
\Sigma_{\bar{M}}$ on $(\bar{M},\Sigma_{\bar{M}})$ is called the
\textit{distribution} of $Z$. If $A$ consists of one element,
i.e.~$A=\{a\}$, I often write $P\{Z=a\}$ instead of $P\{Z\in A\}$.
\begin{definition}\label{stochpro} A $\mathrm{stochastic}$ $\mathrm{process}$
$\{Z_{t};\,t\in\field{Z}\}$ is a one-parameter family of random
variables $Z_{t},\,\,t\in\field{Z}$, defined on the same
probability space $(\Omega,\Sigma_{\Omega},\nu)$ and taking values
in the same measurable space $(\bar{M},\Sigma_{\bar{M}})$.
\end{definition}
The set $\bar{M}$ is called the \textit{outcome space} of the
stochastic process. The bi-infinite sequence $r(\omega):=(\ldots
Z_{-1}(\omega),Z_{0}(\omega),Z_{1}(\omega)\ldots)$ for
$\omega\in\Omega$ is called a \textit{realisation}
\cite[cf.][pp.~4--46]{Doob1953}. Intuitively, $t$ represents time;
so that each $\omega\in\Omega$ represents a possible history in
all its details, and $r(\omega)$ represents the description of that history by giving the score at each $t$.

I will often be concerned with stationary stochastic processes.
These are processes whose probabilistic laws do not change with
time:
\begin{definition} A stochastic process $\{Z_{t};\, t\in \field{Z}\}$ is $\mathrm{stationary}$ if and only if the distributions of the multi-dimensional random variable $(Z_{t_{1}+h},\ldots, Z_{t_{n}+h})$ is the same as the one of $(Z_{t_{1}},\ldots, Z_{t_{n}})$ for all $t_{1}, \ldots, t_{n}\in \field{Z}$, $n\in\field{N}$, and all $h\in \field{Z}$ \cite[p.~94]{Doob1953}.
\end{definition}

It is perhaps needless to stress the importance of discrete
stochastic processes, and stationary processes in particular: both
are ubiquitous in science.

The following stochastic processes will accompany us throughout
the paper. They are probably the most widely known.\\

\noindent\textbf{Example~1: Bernoulli processes}.\\ A Bernoulli
process is a process where, intuitively, at each time point a (possibly biased)
$N$-sided die is tossed where the probability for obtaining side
$s_{k}$ is $p_{k},\,\,1\leq k \leq N,\,\,N\in\field{N}$, with
$\sum_{k=1}^{N}p_{k}\!=\!1$, and each toss is independent of all
the other ones. The mathematical definition proceeds as follows.
The random variables $X_{1},\ldots, X_{n}$, $n\in\field{N}$, are
independent if and only if $P\{X_{1}\in A_{1},\ldots, X_{n}\in
A_{n}\}=P\{X_{1}\in A_{1}\}\ldots P\{X_{n}\in A_{n}\}$ for all
$A_{1},\ldots,A_{n}\in\Sigma_{\bar{M}}$. The random variables
$\{Z_{t};\,t\in\field{Z}\}$ are independent if and only if any
finite number of them is independent.
\begin{definition}\label{Bern}
$\{Z_{t};\,t\in\field{Z}\}$ is a $\mathrm{Bernoulli}$ $\mathrm{process}$ if and only if
(i) its outcome space is a finite number of symbols
$\bar{M}=\{s_{1},\ldots,s_{N}\}, N\in\field{N}$, and
$\Sigma_{\bar{M}}=\field{P}(\bar{M})$, where $\field{P}(\bar{M})$
is the power set of $\bar{M}$; (ii) $P\{Z_{t}=s_{k}\}=p_{k}$ for
all $t\in\field{Z}$ and all $k,1\leq k\leq N$; and (iii)
$\{Z_{t};\,t\in\field{Z}\}$ are independent.
\end{definition}
Clearly, a Bernoulli process is stationary.

In this definition the probability space $\Omega$ is not
explicitly given. I now give a representation of Bernoulli
processes where $\Omega$ is explicitly given. The idea is that
$\Omega$ is the set of realisations of the process. For a Bernoulli process with outcomes $\bar{M}=\{s_{1},\ldots,s_{N}\}$ which have probabilities $p_{1},\ldots,p_{N}$, $N\in\field{N}$, let $\Omega$ be
the set of all sequences $\omega=(\ldots
\omega_{-1}\omega_{0}\omega_{1}\ldots)$ with
$\omega_{i}\in\bar{M}$ corresponding to one of the possible
outcomes of the $i$-th trial in a doubly infinite sequence of
trials. Let $\Sigma_{\Omega}$ be the $\sigma$-algebra generated by
the cylinder-sets
\begin{equation}\label{cylinder}
C^{A_{1}...A_{n}}_{i_{1}...i_{n}}\!\!=\!\!\{\omega\in
\Omega\,|\,\omega_{i_{1}}\!\!\in\!\!A_{1},\!\ldots\!,\omega_{i_{n}}\!\!\in\!\!A_{n},
A_{j}\!\in\!\Sigma_{\bar{M}},\,i_{j}\!\in\!\field{Z},\,i_{1}\!<\!\ldots\!<\!i_{n},\,1\!\leq j\!\leq\!n\}.
\end{equation}
Since the outcomes are independent, these sets have probability
$\bar{\nu}(C^{A_{1}...A_{n}}_{i_{1}...i_{n}}):=P\{Z_{i_{1}}\in
A_{1}\}\ldots P\{Z_{i_{n}}\in A_{n}\}$. Let $\nu$ be defined as the unique extension of
$\bar{\nu}$ to a measure on $\Sigma_{\Omega}$. Finally, define $Z_{t}(\omega):=\omega_{t}$ (the $t$-th
coordinate of $\omega$). Then $\{Z_{t};\,t\in\field{Z}\}$ is the Bernoulli process we started with.

\subsection{Deterministic systems}
According to the canonical definition, a description is
\textit{deterministic} exactly if any two solutions that agree at
one time agree at all times \citep{Butterfield2005}. I
call a sequence which describes the evolution of a deterministic
description over time a solution.

This paper is concerned with measure-theoretic deterministic
descriptions, in short deterministic systems:
\begin{definition}\label{wu}
A $\mathrm{deterministic}$ $\mathrm{system}$ is a quadruple
$(M,\Sigma_{M},\mu,T)$ consisting of a probability space
$(M,\Sigma_{M},\mu)$ and a bijective measurable function
$T$\nolinebreak$:M$\nolinebreak$\rightarrow$\nolinebreak$M$ such that also $T^{-1}$ is measurable.
\end{definition}
The \textit{solution} through $m,\,\,m\in M$, is the sequence
$(T^{t}(m))_{t\in \field{Z}}$. $M$ is the set of all possible
states called the \textit{phase space}; and $T$, which describes
how solutions evolve, is called the \textit{evolution function}.
Clearly, Definition~\ref{wu} defines systems which are
deterministic according to the above canonical definition.

When observing a deterministic system, one observes a value
functionally dependent on, but maybe different from, the actual
state. Hence observations can be
modeled by an \textit{observation function}, i.e.~a measurable
function $\Phi:M\rightarrow M_{O}$ from $(M,\Sigma_{M})$ to the
measurable space $(M_{O},\Sigma_{M_{O}})$
\citep[cf.][p.~16]{OrnsteinWeiss1991}.

I will often be concerned with measure-preserving
deterministic systems \cite[cf.][pp.~3--5]{Cornfeldetal1982}.
\begin{definition}\label{bald}
A $\mathrm{measure}$-$\mathrm{preserving}$ $\mathrm{deterministic}$ $\mathrm{system}$ is a deterministic
system $(M,\Sigma_{M},\mu,T)$ where the measure $\mu$ is
\textit{invariant}, i.e.~for all $A\in\Sigma_{M}$
\begin{equation}\label{invariant}
\mu(T(A))=\mu(A).
\end{equation}
\end{definition}
Measure-preserving deterministic systems are important models in physics but are also important in other sciences such as biology, geology etc. This is so because
condition (\ref{invariant}) is not very restrictive. For first,
all deterministic Hamiltonian systems and statistical-mechanical
systems, and their discrete versions, are measure-preserving; and the relevant
invariant measure is the Lebesgue-measure or a close cousin of it \cite[pp.~5--6]{Petersen1983}.
Second, an invariant measure need not be the Lebesgue measure,
i.e.\ measure-preserving deterministic systems need not be
volume-preserving. Indeed, systems which are not
volume-preserving (called `dissipative systems') can often be
modeled as measure-preserving systems. For instance, the
long-term behaviour of a large class of deterministic systems can be modeled by measure-preserving systems
\citep{EckmannRuelle1985}. More generally, the potential scope of
measure-preserving deterministic systems is quite wide: although
some evolution functions do not have invariant
measures, for very wide classes of evolution functions invariant measures are
proven to exist. For instance, if $T$ is a continuous function on a
compact metric space, there exists at least one invariant measure
\cite[p.~52]{Mane1987}.

I adopt the common assumption that invariant measures can be
interpreted as probability measures. This deep issue has been
discussed in statistical mechanics but is not the focus of this
paper. I only mention two interpretations that naturally suggest
interpreting measures as probability. According to the
time-average interpretation, the measure of a set $A$ is the
long-run average of the time that a solution spends in $A$.
According to the ensemble interpretation, the measure of a set $A$
at $t$ corresponds to the fraction of solutions starting from some
set of initial conditions that are in $A$ at time $t$
[cf.~\citeauthor{EckmannRuelle1985},~\citeyear{EckmannRuelle1985},
pp.~625-627; \citeauthor{Lavis2009}, forthcoming].

The following deterministic system will accompany us.\\

\begin{figure}
\centering
\includegraphics{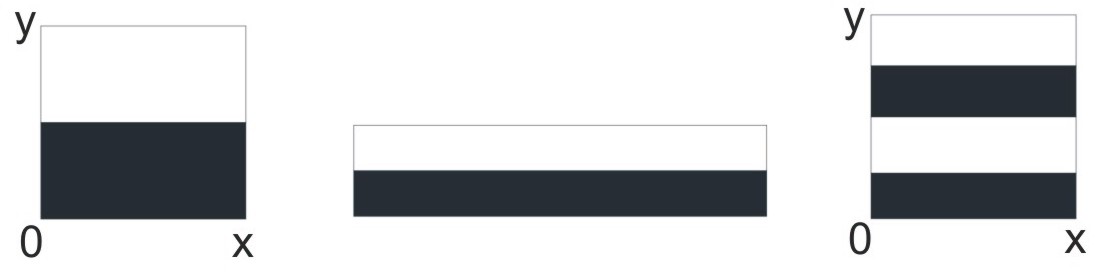}
\caption{\small{The baker's system on $0\leq y \leq
1/2$}\label{Baker}}
\end{figure}

\noindent\textbf{Example~2: The baker's system.}\\
On the unit square $M=[0,1]\times[0,1]\setminus \Gamma$,
where $\Gamma=\{(x,y)\,|\,x=j/2^{n}$ or $y=j/2^{n},
n\in\field{N},0\leq j\leq 2^{n}\}$ consider
\begin{equation}
T(x,y)=(2x,\frac{y}{2})\,\,\,\textnormal{if}\,\,\,0\leq
x<\frac{1}{2};\,\,(2x-1,\frac{y+1}{2})\,\,\,\textnormal{if}\,\,\,\frac{1}{2}\leq
x\leq 1.\end{equation} ($\Gamma$ is excluded to be able to define a bijective function $T$). Figure~1\label{wuzi} illustrates that the
baker's system first stretches the unit square to twice its length
and half its width; then it cuts the rectangle obtained in half
and places the right half on top of the left. For the Lebesgue
measure $\mu$ and the Lebesgue $\sigma$-algebra $\Sigma_{M}$ one
obtains the measure-preserving deterministic system
$(M,\Sigma_{M},\mu,T)$. This system also has physical meaning:
e.g.\ it describes the movement of a particle with initial position $(x,y)$ in the unit square. The particle moves with constant speed perpendicular to the unit square. It bounces on several mirrors, causing it to return to the unit square at $T(x,y)$ \cite[p.~166]{Pitowsky1995}.

\section{Basic observational equivalence}\label{BI}
Let me turn to some results about observational equivalence which
are basic in the sense that they are about the question whether, given a
deterministic system, it is possible to find \textit{any}
stochastic process which is observationally equivalent to the
deterministic system, and conversely.

How can a stochastic process and a deterministic system yield the
same predictions? When a deterministic system is observed, one only
sees how one observed value follows the next observed value. Because the observation function can map two or more actual states
to the same observed value, the same present observed value can
lead to different future observed values. And so a stochastic
process can be observationally equivalent to a deterministic
system only if it is assumed that the deterministic system is
observed with an observation function which is many to one. Yet
this assumption is usually unproblematic: the main reason being that perhaps deterministic systems used in science typically have an infinitely large phase space, and scientists can only observe finitely many different values.

A probability measure is defined on a deterministic system. Hence the predictions derived from a deterministic system are the probability distributions over sequences of possible observations. And similarly, the predictions obtained from a stochastic process are the probability distributions over sequences of possible outcomes. Consequently, the most natural meaning of the phrase \textit{`a stochastic
process and a deterministic system are observationally equivalent'} is:
\textit{(i) the set of possible outcomes of the stochastic process is identical to the the set of possible observed values of the deterministic system, and (ii) the
realisations of the stochastic process and the solutions of the
deterministic system coarse-grained by the observation function
have the same probability distribution.}

Let me now investigate when deterministic systems can be replaced by
stochastic processes. Then I will investigate when stochastic processes
can be replaced by deterministic systems.

\subsection{Deterministic systems replaced by stochastic processes}\label{DRS}

Let $(M,\Sigma_{M},\mu,T)$ be a deterministic system.
According to the canonical Definition~\ref{stochpro},
$Z_{t}(x):=T^{t}(x)$ is a stochastic process with exactly the same
predictions as the deterministic system. However, this process
is evidently equivalent to the original deterministic system, and
the probabilities that one value leads to another one are
trivial ($0$ or $1$). Hence it is still ``really'' a
deterministic system.

But one can do better by appealing to observation functions as
explained above; and, to my knowledge, these results are unknown in
philosophy. Assume the deterministic system
$(M,\Sigma_{M},\mu,T)$ is observed with $\Phi:M\rightarrow M_{O}$.
According to Definition~\ref{stochpro},
$\{Z_{t}:=\Phi(T^{t});\,t\in\field{Z}\}$ is a stochastic process.
It is constructed by applying $\Phi$ to the deterministic system.
Hence the outcomes of the stochastic process are the observed
values of the deterministic system; and the realisations of the
process and the solutions of the deterministic system
coarse-grained by the observation function have the same
probability distribution. Consequently, according to the
characterisation above, \textit{$(M,\Sigma_{M},\mu,T)$ observed
with $\Phi$ is observationally equivalent to stochastic
process $\{\Phi(T^{t});\,t\in\field{Z}\}$}. But the important
question is whether $\{\Phi(T^{t});\,t\in\field{Z}\}$ is
nontrivial. Indeed,  \textit{the stochastic process
$\{\Phi(T^{t});\,t\in\field{Z}\}$ is often nontrivial}. I show now
one result in this direction; besides, several other results also indicate this
\citep[cf.][pp.~178-179]{Cornfeldetal1982}.\footnote{For instance,
if K-systems are observed with a finite-valued observation
function, one obtains nontrivial stochastic processes because for
K-systems the entropy of any finite partition is positive
\citep[cf.][p.~63]{Petersen1983}.}
---
Before I can proceed, the following definitions are needed:
\begin{definition}\label{er}
A measure-preserving deterministic system $(M,\Sigma_{M},\mu,T)$
is $\mathrm{ergodic}$ if and only if for all $A,B\in\Sigma_{M}$
\begin{equation}\label{erg}
\lim_{n\rightarrow\infty}\frac{1}{n}\sum_{i=0}^{n-1}\left(\mu(T^{n}(A)\cap
B)-\mu(A)\mu(B)\right)=0.
\end{equation}
\end{definition}
A \textit{partition} of a measure space $(M,\Sigma_{M},\mu)$ is a
set $\alpha=\{\alpha_{1},\ldots,\alpha_{n}\}$ with
$\alpha_{i}\in\Sigma_{M},\,\,n\in\field{N}$, such that
$\bigcup_{i=1}^{n}\alpha_{i}=M,\,\,\mu(\alpha_{i})>0$, and
$\alpha_{i}\cap\alpha_{j}=\emptyset$ for $i\neq j,\,\,0\leq i,j
\leq n$. A partition is \textit{nontrivial} if and only if it has more than one element.
Let me make the realistic
assumption that the observations have finite accuracy, i.e.~that
only finitely many values are observed. Then one has a
\textit{finite-valued observation function} $\Phi$; i.e.\
$\Phi(m)=\sum_{i=1}^{n}o_{i}\chi_{\alpha_{i}}(m)$,
$M_{O}:=\{o_{i}\,|\,1\leq i\leq n\}$ for some partition $\alpha$
of $(M,\Sigma_{M},\mu)$ and some $n\in\field{N}$, where $\chi_{A}$
denotes the characteristic function of $A$. A finite-valued
observation function is called \textit{nontrivial} if and only if
its corresponding partition is nontrivial
\citep[cf.][p.~179]{Cornfeldetal1982}.

The following proposition shows that for ergodic deterministic
systems for which there is no nontrivial set which is eventually
mapped onto itself, and \textit{every} finite-valued observation
function, the stochastic process $\{\Phi(T^{t});\,t\in\field{Z}\}$
is nontrivial. That is, there is an observed value $o_{i}\in M_{O}$ such that for all observed values $o_{j}\in M_{O}$ the probability of moving from $o_{i}$ to $o_{j}$ is smaller than $1$. Hence there are two or more observed values that can follow  $o_{i}$; and the probability that $o_{i}$ moves to any of these observed values is between $0$ and $1$. This is a strong result because irrespective
of how detailed one looks at the deterministic system, one always
obtains a nontrivial stochastic process.

\begin{proposition}\label{ep}
Assume that the deterministic system $(M,\Sigma_{M},\mu,T)$ is
ergodic and that there does not exist an $n\in\field{N}$ and a
$C\in\Sigma_{M},\,\,0<\mu(C)<1,$ such that, except for
a set of measure zero, $T^{n}(C)=C$. Then for every nontrivial finite-valued
observation function $\Phi:M\rightarrow M_{O}$ and the stochastic
process $\{Z_{t}\!:=\!\Phi(T^{t});\,t\in\field{Z}\}$ the following
holds: there is an $o_{i}\in M_{O}$ such that for all
$o_{j}\in M_{O}$, $P\{Z_{t+1}\!=\!o_{j}\,|\,Z_{t}\!=\!o_{i}\}<1$.\footnote{For a
random variable $Z$ to a measurable space
$(\bar{M},\Sigma_{\bar{M}})$ where $\bar{M}$ is finite the
conditional probability is defined as usual as:\\ $P\{Z\in
A\,|\,Z\in B\}:=P\{Z\in A\cap B\}/P\{Z\in B\}$ for all $A,B\subseteq
\Sigma_{\bar{M}}$ with $P\{Z\in B\}>0$.}
\end{proposition}
For a proof, see subsection~\ref{PP1}. For instance, the baker's system (Example~2) is weakly mixing, and thus any finite-valued observation function gives rise to a nontrivial stochastic process.

Measure-preserving systems are typically what is called `weakly mixing'\footnote{
$(M,\Sigma_{M},\mu,T)$ is weakly mixing if and only if for all
$A,B\in\Sigma_{M}$
\begin{equation}\label{wm}
\lim_{n\rightarrow\infty}\frac{1}{n}\sum_{i=0}^{n-1}|\mu(T^{n}(A)\cap
B)-\mu(A)\mu(B)|=0.
\end{equation}}
\citep{Halmos1944}. It is easy to see that any weakly mixing
system satisfies the assumption of Proposition~\ref{ep} (weakly mixing is stronger than this assumption).\footnote{First, it is clear that weakly mixing systems
are ergodic. Second, assume that for a weakly mixing system there
exists an $n\in\field{N}$ and a $C\in\Sigma_{M},$ $0<\mu(C)<1,$
such that, except for a set of measure zero, $T^{n}(C)=C$. But then
equation~(\ref{wm}) cannot hold for $A:=C$ and $B:=C$. In
subsection~\ref{PM2} I will show that the irrational rotation
on the circle satisfies the assumption of
Proposition~\ref{ep} but is not weakly mixing.}
Hence Proposition~\ref{ep} shows that for typical measure-preserving
deterministic systems any finite-valued observation function
yields a nontrivial stochastic process.

Yet Proposition~1 does not say much about whether the measure-preserving deterministic systems encountered in science fulfill the assumption of Proposition~\ref{ep} because those systems constitute a small class of all measure-preserving systems. Indeed, the KAM theorem says that the phase space of integrable Hamiltonian systems which are perturbed by a small nonintegrable perturbation breaks up into stable and unstable regions. With increasing perturbation the unstable regions become larger and often eventually cover nearly the entire hypersurface of constant energy. Because a solution is confined to a region, the KAM theorem implies that the discrete versions of infinitely differentiable Hamiltonian systems are typically not ergodic \citep[section 4]{Berkovitzetal2006}. (I call the discrete-time systems obtained by looking at a continuous-time system $S$ at points of time $nt_{0}$, $n\in\field{N}$, $t_{0} \in \field{R}$ arbitrary, $t_{0}\neq 0$, the \textit{discrete versions of $S$}).\footnote{Alternatively, continuous-time deterministic systems can be discretised by considering the successive hits of a trajectory on a suitable Poincar\'{e} section. All I say about discrete versions of continuous systems also holds true for discrete-time systems arising in this way, except that the latter are more often ergodic \cite[pp.~680--685]{Berkovitzetal2006}.}

Despite this, Proposition~\ref{ep} applies to several systems encountered in science.
First, a motion is chaotic just in case it is deterministic yet also unstable because nearby initial conditions eventually lead to very different outcomes. I will not need a more exact definition; but I will call a system \textit{chaotic} if the motion is chaotic on the entire phase space and \textit{locally chaotic} if the motion is chaotic on a region of phase space. Chaotic systems are usually regarded as weakly mixing (\citeauthor{Berkovitzetal2006}, \citeyear{Berkovitzetal2006}, p.~688; \citeauthor{Werndl2009a}, \citeyear{Werndl2009a}, section 3). And as will be argued in subsection~\ref{DSRSS}, there are several physically relevant chaotic and weakly mixing systems. Moreover, in subsection~\ref{PM2} it will be shown that there are even systems which are neither chaotic nor locally chaotic but which satisfy Proposition~\ref{ep}.
Second, even if the whole system does not satisfy the assumption of Proposition~\ref{ep}, the motion of the system restricted to some regions of phase space might well satisfy this assumption. In fact, Proposition~\ref{ep} immediately implies the following result. Assume that for a measure-preserving system $(M,\Sigma_{M},\mu,T)$ there is a $A\in\Sigma_{M}, \mu(A)>0$, such that the system restricted to $A$\footnote{That is, the system $(A,\Sigma_{M\cap A},\mu_{A},T_{A})$, where
$\Sigma_{M\cap A}:=\{B\cap A|\,B\in\Sigma_{M}\},$ $\mu_{A}(X):=\frac{\mu(X)}{\mu(A)}$, and $T_{A}$ denotes $T$ restricted to $A$.} fulfills the assumption of Proposition~\ref{ep}. Then all observations which discriminate between values in $A$ lead to nontrivial stochastic processes. That is, for any observation function $\Phi(m)=\sum_{i=1}^{n}o_{i}\chi_{\alpha_{i}}(m)$ where there are $k,l,\, k\neq l$, such that $\mu(A\cap\alpha_{k})\neq 0$ and $\mu(A\cap\alpha_{l})\neq 0$, there is an outcome $o_{i}\in M_{O}$ such that for all outcomes $o_{j}\in M_{O}$ it holds that $P\{Z_{t+1}\!=\!o_{j}\,|\,Z_{t}\!=\!o_{i}\}<1$.
In particular, although mathematically little is known, it is conjectured that the motion restricted to unstable regions of KAM-type systems is weakly mixing (\citeauthor{Berkovitzetal2006}, \citeyear{Berkovitzetal2006}, section 4; \citeauthor{Werndl2009a}, \citeyear{Werndl2009a}, section 3). If this is true, then my argument shows that for many observation functions of KAM-type systems one obtains nontrivial stochastic processes.

\subsection{Stochastic processes replaced by deterministic systems}\label{SRD}

I have shown that deterministic systems, when observed, can yield
nontrivial stochastic processes. But can one find, for every stochastic
process, a deterministic system which produces this process?

The following idea of how to replace stochastic processes by
deterministic systems is well known in the technical literature
\citep[pp.~6--7]{Petersen1983}\footnote{\citeauthor{Petersen1983}
discusses it only for stationary stochastic processes; I consider
generally stochastic processes.} and known to philosophers
\citep{Butterfield2005}; I also need to discuss it for what
follows later. The underlying thought is that for each realisation $r(\omega)$, one sets up a deterministic
system with phase space $\{r(\omega)\}$.

So consider a stochastic
process $\{Z_{t};\,t\in\field{Z}\}$ from
$(\Omega,\Sigma_{\Omega},\nu)$ to $(\bar{M},\Sigma_{\bar{M}})$.
Let $M$ be the set of all bi-infinite sequences $m=(\ldots
m_{-1}m_{0}m_{1}\ldots)$ with $m_{i}\in\bar{M},\,\,i\in\field{Z}$,
and let $m_{t}$ be the $t$-th coordinate of $m,\,\,t\in\field{Z}$.
Let $\Sigma_{M}$ be the $\sigma$-algebra generated by the cylinder-sets as defined in $(\ref{cylinder})$ at the end of subsection~\ref{SP}.
$\{Z_{t}; t\in\field{Z}\}$ assigns to each cylinder set $C^{A_{1}...A_{n}}_{i_{1}...i_{n}}$ a pre-measure,
namely the probability $P\{Z_{i_{1}}\in
A_{1},\ldots,Z_{i_{n}}\in A_{n}\}$. Let $\mu$ be the unique
extension of this pre-measure to a measure on $\Sigma_{M}$. Let
$T:M\rightarrow M$ be the left shift, i.e.~$T((\ldots
m_{-1}m_{0}m_{1}\ldots)):=(\ldots m_{0}m_{1}m_{2}\ldots)$. $T$ is
bijective and measurable, and so one obtains the deterministic
system $(M,\Sigma_{M},\mu,T)$. Finally, assume one sees only the
0-th coordinate of the sequence $m$, i.e.\ one applies the
observation function
$\Phi_{0}:M\rightarrow\bar{M},\Phi_{0}(m)=m_{0}$. I now define:
\begin{definition}\label{DetRep}
$(M,\Sigma_{M},\mu,T,\Phi_{0})$ as constructed above is the
$\mathrm{deterministic}$ \linebreak $\mathrm{representation}$ of the process
$\{Z_{t};\,t\in\field{Z}\}$.
\end{definition}

For the deterministic representation
$(M,\Sigma_{M},\mu,T,\Phi_{0})$ of $\{Z_{t}; t\in\field{Z}\}$ it is assumed that the 0-th coordinate is observed. Consequently, the possible outcomes
of $\{Z_{t}; t\in\field{Z}\}$ are the possible observed values of
$(M,\Sigma_{M},\mu,T,\Phi_{0})$. Clearly, any realisation $r(\omega)$ of the process is contained in $M$, and observing the solution $(T^{t}(r(\omega)))_{t\in\field{Z}}$
with $\Phi_{0}$ exactly gives $r(\omega)$.
Furthermore, the measure $\mu$ is defined by the probabilities
which are assigned by $\{Z_{t}; t\in\field{Z}\}$ to each cylinder
set. Hence the probability distribution over the realisations of
$\{Z_{t}; t\in\field{Z}\}$ is the same as the one over the sequences of observed values of $(M,\Sigma_{M},\mu,T,\Phi_{0})$. Thus, according to the characterisation at the start of this section, \textit{a
stochastic process is observationally equivalent to its
deterministic representation. Hence every stochastic process can
be replaced by at least one deterministic system.} (When there is
no risk of confusion, I also refer to the system
$(M,\Sigma_{M},\mu,T)$ of the deterministic representation
$(M,\Sigma_{M},\mu,T,\Phi_{0})$ as the deterministic
representation.)

For Bernoulli processes (Example~1) the deterministic
representation is the following. $(M,\Sigma_{M},\mu)$ is the
measure space $(\Omega,\Sigma_{\Omega},\nu)$ as defined
 at the end of subsection~\ref{SP}. $T((\ldots\omega_{-1}\omega_{0}\omega_{1}\ldots)):=(\ldots\omega_{0}\omega_{1}\omega_{2}\ldots)$
for $\omega\in\Omega$ and $\Phi_{0}(\omega)=\omega_{0}$.
\begin{definition}\label{BS}
The deterministic representation $(M,\Sigma_{M},\mu,T)$ of the
Bernoulli process with probabilities
$p_{1},\ldots,p_{N},\,\,N\in\field{N},$ is called the
$\mathrm{Bernoulli}$ $\mathrm{shift}$ with probabilities
$(p_{1},\ldots,p_{N})$.
\end{definition}

From a philosophical perspective the deterministic representation
is a cheat because its states are constructed to encode the future
and past outcomes of the stochastic process. Despite this, it is important to know that the deterministic
representation exists. Of course, there is the question whether
deterministic systems which do not involve a cheat can replace a stochastic process. I will turn to this question in section~\ref{AI}, where I show that for some stochastic processes this is indeed the case. To my knowledge, it is unknown whether \textit{every} stochastic process can be thus replaced.

\subsection{A mathematical definition of observational equivalence}
Let me now mathematically define what it means for a
stochastic process and a deterministic system to be
observationally equivalent. The notion of isomorphism captures the idea that deterministic systems are probabilistically equivalent, i.e.~that their states can be put into
one-to-one correspondence such that the corresponding solutions
have the same probability distributions.
\begin{definition}\label{isomorphic}
$(M_{1},\Sigma_{M_{1}},\mu_{1},T_{1})$ is $\mathrm{isomorphic}$ to
$(M_{2},\Sigma_{M_{2}},\mu_{2},T_{2})$ (where both systems are
assumed to be measure-preserving) if and only if there are
measurable sets $\hat{M}_{i}\subseteq M_{i}$ with
$\mu_{i}(M_{i}\setminus \hat{M}_{i})=0$ and
$T_{i}\hat{M}_{i}\subseteq\hat{M}_{i}\,\,(i=1,2$), and there is a
bijection $\phi:\hat{M}_{1}\!\rightarrow \!\hat{M}_{2}$ such that
(i) $\phi(A)\!\in\!\Sigma_{M_{2}}$ for all
$A\!\in\!\Sigma_{M_{1}},A\subseteq \hat{M}_{1}$, and
$\phi^{-1}(B)\in\Sigma_{M_{1}}$ for all
$B\in\Sigma_{M_{2}},B\subseteq \hat{M}_{2}$; (ii)
$\mu_{2}(\phi(A))=\mu_{1}(A)$ for all
$A\in\Sigma_{M_{1}},\,A\subseteq\hat{M}_{1}$; (iii)
$\phi(T_{1}(m))=T_{2}(\phi(m))$ for all $m\in\hat{M}_{1}$
\citep[cf.][p.~4]{Petersen1983}.
\end{definition}
One easily sees that `being isomorphic' is an equivalence
relation. Isomorphic systems may have different phase spaces. If
identical sets $\hat{M}_{1}$ and $\hat{M}_{2}$ can be found, then
the deterministic systems are obviously probabilistically equivalent and
have, from a probabilistic viewpoint, the same phase space; for this case it will later be convenient to say that the measure-preserving deterministic systems are
\textit{manifestly isomorphic.}

According to the characterisation at the beginning of this
section, a deterministic system $(M,\Sigma_{M},\mu,T)$, observed
with $\Phi$, gives the same predictions as
$\{Z_{t}|\,t\in\field{Z}\}$ exactly if (i) the outcomes of
$\{Z_{t}|\,t\in\field{Z}\}$ are the observed values of
$(M,\Sigma_{M},\mu,T)$, and (ii) the deterministic representation
of $\{\Phi(T^{t});\,t\in \field{Z}\}$ is probabilistically
equivalent to the deterministic representation of
$\{Z_{t}|\,t\in\field{Z}\}$. Hence one arrives
at the following definition of `observational equivalence'; (for
what follows, a definition for measure-preserving systems and,
correspondingly, stationary stochastic processes will
suffice):\footnote{For a measure-preserving system
$(M,\Sigma_{M},\mu,T)$ the process
$\{\Phi(T^{t});\,t\in\field{Z}\}$ is stationary: $\{x\in
M\,|\,\Phi(T^{t_{1}}(x))\in A_{1},\ldots,\Phi(T^{t_{n}}(x))\in
A_{n},\,A_{i}\in
\Sigma_{M_{O}},t_{i}\in\field{Z},\,n\in\field{N}\}$ is identical
to $A:=T^{-t_{1}}(\Phi^{-1}(A_{1})\cap\ldots\cap
T^{t_{1}-t_{n}}\Phi^{-1}(A_{n}))$. Likewise, $\{x\in
M\,|\,\Phi(T^{t_{1}+h}(x))\in A_{1},\ldots,\Phi(T^{t_{n}+h}(x))\in
A_{n}\}$ is $B:=T^{-(t_{1}+h)}(\Phi^{-1}(A_{1})\cap\ldots\cap
T^{t_{1}-t_{n}}\Phi^{-1}(A_{n}))$. Because the system is
measure-preserving, $\mu(A)=\mu(B)$, implying that
$\{\Phi(T^{t});\,t\in \field{Z}\}$ is stationary. And if
$\{Z_{t}|\,t\in\field{Z}\}$ is stationary, its deterministic
representation $(M,\Sigma_{M},\mu,T)$ is measure-preserving. For
stationarity implies that $\mu(T(A))=\mu(A)$ for any cylinder set
$A$, and hence that $\mu(T(A))=\mu(A)$ for all $A\in\Sigma_{M}$
\citep[cf.][p.~178]{Cornfeldetal1982}.}
\begin{definition}\label{replace1}
The stationary stochastic process $\{Z_{t};\,t\in\field{Z}\}$ and
the measure-preserving deterministic system
$(M,\Sigma_{M},\mu,T)$, observed with $\Phi$, are $\mathrm{observationally}$
$\mathrm{equivalent}$ if and only if the deterministic representation of
$\{\Phi(T^{t});\,t\in \field{Z}\}$ is manifestly isomorphic to the
deterministic representation of $\{Z_{t};\,t\in\field{Z}\}$.
\end{definition}

All the cases of observational equivalence already discussed are cases
of observational equivalence in the sense of
Definition~\ref{replace1}. First, I claimed in
subsection~\ref{DRS} that $(M,\Sigma_{M},\mu,T)$ observed with
$\Phi$ is observationally equivalent to the stochastic process
$\{\Phi(T^{t});\,t\in\field{Z}\}$. This is true because every system is manifestly isomorphic to itself. Second, I claimed in
subsection~\ref{SRD} that the deterministic representation
$(M,\Sigma_{M},\mu,T,\Phi_{0})$ of $\{Z_{t};\,t\in\field{Z}\}$ is
observationally equivalent to $\{Z_{t};\,t\in\field{Z}\}$. This is
true because the deterministic representation of
$\{\Phi_{0}(T^{t});\,t\in\field{Z}\}$ is
$(M,\Sigma_{M},\mu,T,\Phi_{0})$.

One final point: assume that
$(M,\Sigma_{M},\mu,T)$ is isomorphic via
$\phi:\hat{M}\rightarrow\hat{M_{2}}$ to the deterministic
representation $(M_{2},\Sigma_{M_{2}},\mu_{2},T_{2},\Phi_{0})$ of
$\{Z_{t};\,t\in\field{Z}\}$. This means that \textit{there is a
one-to-one correspondence between the solutions of the
deterministic system and the realisations of the stochastic
process.} Thus $(M,\Sigma_{M},\mu,T)$ restricted to
$\hat{M}$ and observed with $\Phi_{0}(\phi(m))$ is
observationally equivalent to $\{Z_{t};\,t\in\field{Z}\}$. This is
so because the deterministic representation of
$\{\Phi_{0}(\phi(T^{t}));\,t\in\field{Z}\}$ where $T$ is
restricted to $\hat{M}$ is identical to the deterministic
representation of $\{\Phi_{0}(T_{2}^{t});\,t\in\field{Z}\}$ where
$T_{2}$ is restricted to $\hat{M_{2}}$. Hence the deterministic
representation of $\{\Phi_{0}(\phi(T^{t}));\,t\in \field{Z}\}$
is manifestly isomorphic to $(M_{2},\Sigma_{M_{2}},\mu_{2},T_{2})$.

The following definition will be important later:
\begin{definition}\label{BerSys}
$(M,\Sigma_{M},\mu,T)$ is a $\mathrm{Bernoulli}$ $\mathrm{system}$ if and only if
it is isomorphic to a Bernoulli shift.
\end{definition}
The meaning of Bernoulli systems is
clear, viz.\ the solutions of a Bernoulli system can put into
one-to-one correspondence with the realisations of a Bernoulli
process. Thus a Bernoulli system, observed with $\Phi_{0}(\phi)$,
produces a Bernoulli process. Finally, I note the important result that two Bernoulli shifts (and hence two Bernoulli systems) are isomorphic if and only if they have the same Kolmogorov-Sinai entropy, where the Kolmogorov-Sinai entropy of a Bernoulli shift with probabilities $(p_{1},\ldots,p_{n})$ is $\sum_{i=1}^{n}-p_{i}\log_{2} p_{i}$ (\citeauthor{FriggWerndl}, 2010; \citeauthor{Ornstein1974}, \citeyear{Ornstein1974}, pp.~5--3; \citeauthor{Werndl2009b}, \citeyear{Werndl2009b}).

\section{Advanced observational equivalence}\label{AI}

In this section I discuss results which are `advanced' in the sense that they are about the question whether it is possible to replace deterministic systems \textit{in
science} with stochastic processes \textit{in science}. The phrase
`systems in science' (or `processes in science') is a short-hand
for systems (or processes) which are used in science to model phenomena.

\subsection{Deterministic system in science which replace stochastic processes in science}\label{DSRSS}
The deterministic representation does not naturally arise in science (no doubt reflecting that fact that is a philosophical cheat). And the results so far only show that stochastic processes in science,
e.g.\ a Bernoulli process, can be replaced by its deterministic
representation. Hence it seems hard to imagine how deterministic systems in science could replace stochastic processes in science. In particular, it seems hard to
imagine how deterministic systems in science could be random
enough to replace random stochastic processes such as Bernoulli
processes. Thus one might conjecture that \textit{it is impossible to replace stochastic processes in science by deterministic systems in science.}

Bernoulli processes (Example~1) are often regarded as the most
random discrete-time stochastic processes because their outcomes are independent (cf.~\citeauthor{Ornstein1989},
1989). Are there deterministic systems in science
which, when observed, are observationally equivalent to Bernoulli
processes? And are there even deterministic systems in science
which are Bernoulli systems? Historically, it was long thought
that the answer to these questions is negative \citep[cf.][p.~834]{Sinai1989}.
So it was a big surprise when it was discovered from the 1960s onwards that
\textit{there are several deterministic systems in science which
are Bernoulli systems} (among them systems producing Bernoulli
processes with equiprobable outcomes).  Let me mention some of the most important examples, which are also some of the most important examples of chaotic systems.\footnote{Bernoulli systems are regarded as strongly chaotic.}

To start with, there are systems in Newtonian mechanics, some of which are simple models of statistical mechanical systems, whose discrete versions are proven to be Bernoulli systems.  The most prominent examples are: first, some hard sphere systems, which describe the motion of a number of hard spheres undergoing elastic reflections at the boundary and collisions amongst each other; e.g., the motion of $N$ hard balls on the $m$ torus for $N\geq 2$ and $m\geq N$; second, billiard systems with convex obstacles; and third, geodesic flows of negative curvature, i.e.\ frictionless motion of a particle moving with unit speed on a compact manifold with everywhere negative curvature. It is usually very hard to prove that systems are Bernoulli. Therefore, for many systems it is only conjectured that their discrete versions are Bernoulli, e.g., for \textit{all} hard sphere systems and the motion of KAM-type systems restricted to some regions of phase space (\citeauthor{OrnsteinWeiss1991}, \citeyear{OrnsteinWeiss1991}, section 4; \citeauthor{Young1997}, \citeyear{Young1997}; \citeauthor{Berkovitzetal2006}, \citeyear{Berkovitzetal2006}, p.~679--680).

Furthermore, there are dissipative systems which are Bernoulli systems: such as the logistic map and generalised versions thereof, the H\'{e}non map and generalised versions thereof, and the discrete versions of the Lorenz system and generalised versions thereof. Some of these systems give relatively accurate predictions, e.g.\ the Lorenz system as a model for water-wheels. Yet often these systems are motivated as simple models which help us to understand, and not so much to predict, phenomena: e.g.\ the logistic map for population and climate dynamics, and the H\'{e}non map and the Lorenz system for weather dynamics (\citeauthor{Lorenz1964}, \citeyear{Lorenz1964}; \citeauthor{May1976}, \citeyear{May1976}; \citeauthor{Jacobson1981}, \citeyear{Jacobson1981}; \citeauthor{BenedicksYoung1993},~\citeyear{BenedicksYoung1993};
\citeauthor{Smith1998},~\citeyear{Smith1998}, chapter 8;
\citeauthor{Lyubich2002},~\citeyear{Lyubich2002};
\citeauthor{Luzzatto2005}, \citeyear{Luzzatto2005}).

Also the baker's system $(M,\Sigma_{M},\mu,T)$
(Example~2), a somewhat artificial example of deterministic
motion, is Bernoulli. Assign to each $(x,y)$ in $M$
the sequence
$\phi(x,y)=\ldots\omega_{-2}\omega_{-1}\omega_{0}\omega_{1}\omega_{2}\ldots$ defined by the binary expansion of the coordinates:
\begin{equation}
x=0.\omega_{0}\omega_{1}\ldots=\sum_{i=1}^{\infty}\frac{\omega_{i-1}}{2^{i}};\,\,\,\,\,y=0.\omega_{-1}\omega_{-2}\ldots=\sum_{i=1}^{\infty}\frac{\omega_{-i}}{2^{i}}. \end{equation}
Consider the Bernoulli shift
$(M_{2},\Sigma_{M_{2}},\mu_{2},T_{2})$ with states $s_{1},s_{2}$ and probabilities
$(\frac{1}{2},\frac{1}{2})$. Let $\hat{M}_{2}$ be the
subset of $M_{2}$ excluding all states beginning or ending with an infinite sequence of ones or zeros; note that $\mu_{2}(\hat{M}_{2})=1$. One easily
verifies that $\phi:M\rightarrow\hat{M}_{2}$ gives an isomorphism
from $(M,\Sigma_{M},\mu,T)$ to
$(M_{2},\Sigma_{M_{2}},\mu_{2},T_{2})$. Hence the
baker's system with the observation function
$\Phi((x,y)):=s_{1}\chi_{\alpha_{1}}((x,y))+s_{2}\chi_{\alpha_{2}}((x,y))$,
where
$\alpha=\{\alpha_{1},\alpha_{2}\}:=\{[0,\frac{1}{2})\times[0,1]\setminus \Gamma,[\frac{1}{2},1]\times[0,1]\setminus \Gamma\}$ yields the Bernoulli process
with states $s_{1},s_{2}$ and probabilities $(\frac{1}{2},\frac{1}{2})$.

A Bernoulli system is weakly mixing \citep[p.~58]{Petersen1983}.
Hence, provided it is observed with a
finite-valued observation function, one always obtains a
nontrivial stochastic process (Proposition~1).

What is the significance of the these results? They show that the conjecture advanced at the beginning of this subsection is wrong: \textit{it is possible to replace stochastic processes in science by  deterministic systems in science.}\footnote{The arguments in this section allow any meaning of `deterministic systems in science' that is wide enough to include some Bernoulli systems but narrow enough to exclude systems such as the deterministic representation.}

Of course, the question arises whether for deterministic systems
in science which are observationally equivalent to stochastic
processes in science the corresponding observation function is
\textit{natural} in the sense that one might encounter it when modeling
phenomena. The answer depends on the deterministic system and the
phenomenon under consideration. For some systems the observation
function is very involved and thus no natural interpretation can
be found. But in other cases the observation function corresponds
to a realistic way of observing the system.

For instance, recall that the baker's system models a particle bouncing on several
mirrors where $(x,y)$ denotes the position of the particle on a
square. Here an observer might well only be interested in whether
the position of the particle is to the left or to the right of the
square. Then the observation function  $\Phi((x,y)):=s_{1}\chi_{\alpha_{1}}((x,y))+s_{2}\chi_{\alpha_{2}}((x,y))$,
above, which indeed produces a Bernoulli process, would be natural.

\subsection{Stochastic processes which replace deterministic systems in science at every observation level}\label{VC}

\subsubsection{$\varepsilon$-congruence and replacement by Markov processes}\label{muede}

The previous discussion showed that for several deterministic systems in science, regardless which finite-valued observation function one applies, one always obtains  a stochastic process. But to obtain systems in science such as Bernoulli processes, it seems crucial that \textit{coarse} observation functions are applied.
Hence it is hard to imagine that by taking finer and finer
observations of deterministic systems in science one still obtains
stochastic processes in science. In particular, it is hard to
imagine that one still obtains random stochastic processes. Therefore, one might
conjecture that \textit{it is impossible to replace deterministic systems in science at every observation level by stochastic processes in science.}

Let me introduce one of the most natural ways of understanding the
phrase `at any observation level', i.e.\ the notion that stochastic processes of a certain type replace a deterministic system at any observation level. I first explain
what it means for a deterministic system and a stochastic process
to give the same predictions at an observation level $\varepsilon>0$, $\varepsilon \in\field{R}$. There are two aspects. First, one imagines that in practice, for sufficiently small $\varepsilon_{1}$, one cannot
distinguish states of the deterministic system which are less than
the distance $\varepsilon_{1}$ apart. The second aspect concerns probabilities: in practice, for sufficiently small $\varepsilon_{2}$, one will not be able to observe
differences in probabilities of less than
$\varepsilon_{2}$. Assume that $\varepsilon$ is smaller than
$\varepsilon_{1}$ and $\varepsilon_{2}$. Then a deterministic
system and a stochastic process give the same predictions at
observation level $\varepsilon$ if the following holds: the
solutions of the deterministic system can be put into one-to-one
correspondence with the realisations of the stochastic process in such a way that the actual state of the deterministic system and the
corresponding outcome of the stochastic process are at each time point less then
$\varepsilon$ apart except for a set whose probability is
smaller than $\varepsilon$.

Mathematically, this idea is captured by the notion of
$\varepsilon$-congruence. To define it, one needs to speak of
distances between states in the phase space $M$ of the
deterministic system; hence one assumes a metric
$d_{M}$ defined on $M$. So we need to find a
stochastic process whose outcome is within distance $\varepsilon$
of the actual state of the deterministic system. Hence one assumes
that the possible outcomes of the stochastic process are a subset
of the phase space of the deterministic system. For what follows,
it suffices to consider measure-preserving deterministic systems
and, correspondingly, stationary processes. Now recall
Definition~\ref{DetRep} of the deterministic representation and
Definition~\ref{isomorphic} of being isomorphic. So finally, I can define:
\begin{definition}
Let $(M,\Sigma_{M},\mu,T)$ be a measure-preserving
deterministic system, where $(M,d_{M})$ is a metric space.
Let $(M_{2},\Sigma_{M_{2}},\mu_{2},T_{2},\Phi_{0})$ be the
deterministic representation of the stationary stochastic process
$\{Z_{t};\,t\in\field{Z}\}$, which takes values in
$(M,d_{M})$, i.e.~$\Phi_{0}:M_{2}\rightarrow M$.
$(M,\Sigma_{M},\mu,T)$ is
$\varepsilon$-$\mathrm{congruent}$ to
$(M_{2},\Sigma_{M_{2}},\mu_{2},T_{2},\Phi_{0})$ if and only if
$(M,\Sigma_{M},\mu,T)$ is isomorphic via a
function $\phi:M\rightarrow M_{2}$ to
$(M_{2},\Sigma_{M_{2}},\mu_{2},T_{2})$ and
$d_{M}(m,\Phi_{0}(\phi(m)))<\varepsilon$ for all $m\in M$
except for a set of measure $<\varepsilon$ in $M$
\citep[cf.][pp.~22--23]{OrnsteinWeiss1991}.
\end{definition}
Note that $\varepsilon$-congruence does not assume that the deterministic system is observed with an observation function. Of course, observation functions can be introduced. Assume one observes a deterministic system with an observation function. Then  there is a stochastic process which is $\varepsilon$-congruent to the deterministic system such that the probabilistic predictions resulting from the observation function differ at most by $\varepsilon$ from the probabilistic predictions obtained by applying the observation function to the $\varepsilon$-congruent stochastic process.

By generalising over $\varepsilon$, one obtains a natural meaning
of the notion that stochastic processes of a certain type replace
a measure-preserving deterministic system at any observation
level, namely: for every $\varepsilon>0$ there is a stochastic
process of this type which gives the same predictions at
observation level $\varepsilon$. Or technically: \textit{for every
$\varepsilon>0$ there exists a stochastic process of this type
which is $\varepsilon$-congruent to the deterministic system.}

For Bernoulli processes the next outcome of the process is
independent of its previous outcomes. So, intuitively, it seems
clear that deterministic systems in science, for which the next
state of the system is constrained by its previous states (because
of the underlying determinism at the level of states), cannot be
replaced by Bernoulli processes at \textit{every} observation level. \citet[pp.~160--162]{Smith1998} also hints at this idea but does not substantiate
it with a proof. The following theorem shows that,
for our notion of replacement at every observation level, this
idea is indeed correct under very mild assumptions, which hold for deterministic systems in science.\footnote{This theorem also holds for \textit{generalised
Bernoulli processes}---stochastic processes consisting of
independent and identically distributed random variables whose
outcome space need not be finite. That is, under
the assumption of Theorem~\ref{T1}, there is an $\varepsilon>0$
for which there is no generalised Bernoulli process to which
$(M,\Sigma_{M},\mu,T)$ is $\varepsilon$-congruent (cf.~Remark~1
at the end of subsection~\ref{PT1}).} Hence this theorem shows a limitation on the observational equivalence of deterministic systems and stochastic processes.

\begin{theorem}\label{T1}
Let $(M,\Sigma_{M},\mu,T)$ be a measure-preserving deterministic
system where $\Sigma_{M}$ contains all open balls of the metric
space $(M,d_{M})$, $T$ is continuous at some point $x\in M$, every
open ball around $x$ has positive measure, and there is a set
$D\in\Sigma_{M},\,\,\mu(D)>0$, with
$d(T(x),D):=\inf\{d(T(x),m)\,|\,m \in D\}>0$. Then there is some
$\varepsilon>0$ for which there is no Bernoulli process to which
$(M,\Sigma_{M},\mu,T)$ is $\varepsilon$-congruent.
\end{theorem}
For a proof, see subsection~\ref{PT1}.\footnote{A deterministic
system which is replaced at every observation level by a Bernoulli
process will be measure-preserving. Hence it suffices to
concentrate on measure-preserving systems.}

Given this result, it is natural to ask (which, incidentally, \citeauthor{Smith1998} (\citeyear{Smith1998}) does not do) whether deterministic systems in science can be replaced at every observation level by other stochastic processes in science. \textit{The answer is `yes'. Besides, all one needs are
irreducible and aperiodic Markov processes, which are widely used
in science.} These Markov processes are often regarded as random;
in particular, Bernoulli processes are regarded as the most random processes and
Markov processes as the next most random
(\citeauthor{OrnsteinWeiss1991},~\citeyear{OrnsteinWeiss1991},
p.~38 and p.~66).

For Markov processes the next outcome depends only on the previous
outcome.
\begin{definition} $\{Z_{t};\,t\in\field{Z}\}$ is a
$\mathrm{Markov}$ $\mathrm{process}$ if and only if (i) its outcome space consists of a
finite number of symbols
$\bar{M}:=\{s_{1},\ldots,s_{N}\},\,\,N\in\field{N}$, and
$\Sigma_{\bar{M}}=\field{P}(\bar{M})$; (ii)
$P\{Z_{t+1}=s_{j}\,|\,Z_{t},Z_{t-1}\ldots,Z_{k}\}=P\{Z_{t+1}=s_{j}\,|\,Z_{t}\}$
for any $t$, any $k\in\field{Z},\,k\leq t$, and any
$s_{j}\in\bar{M}$; and (iii) $\{Z_{t};\,\,t\in\field{Z}\}$ is stationary
\end{definition}

Define $P^{k}(s_{i},s_{j}):=P\{Z_{n+k}=s_{i}\,|\,Z_{n}=s_{j}\}$
for $k\in\field{Z}$. A Markov process is \textit{irreducible}
exactly if it cannot be split into two processes because each
outcome can be reached from all other outcomes;
formally: for every $s_{i},s_{j}\in\bar{M}$ there is a
$k\in\field{N}$ such that $P^{k}(s_{i},s_{j})>0$. A Markov process
is aperiodic exactly if for every possible outcome there is no
periodic pattern in which the process can visit that outcome.
Mathematically, the \textit{period} $d_{s_{i}}$ of an outcome
$s_{i}\in\bar{M},\,\,1\leq i\leq N$, is defined by
$d_{i}={gcd}\{n\geq 1\,|\,P^{n}(s_{i},s_{i})>0\}$ where `$gcd$'
denotes the greatest common divisor. An outcome $s_{i}\in\bar{M}$
is \textit{aperiodic} if and only if $d_{i}=1$, and the Markov
process is \textit{aperiodic} if and only if all its possible
outcomes are aperiodic.

The following theorem shows that Bernoulli systems (cf.~Definition
\ref{BerSys}) can be replaced at every observation level by
irreducible and aperiodic Markov processes.
\begin{theorem}\label{T2}
Let $(M,\Sigma_{M},\mu,T)$ be a Bernoulli system where the metric
space $(M,d_{M})$ is separable\footnote{$(M,d_{M})$ is separable
if and only if there exists a countable set
$\ddot{M}=\{m_{n}\,|n\in\field{N}\}$ with $m_{n}\in M$ such that
every nonempty open subset of $M$ contains at least one element of
$\ddot{M}$.} and $\Sigma_{M}$ contains all open balls of
$(M,d_{M})$. Then for any $\varepsilon>0$ there is an irreducible
and aperiodic Markov process such that $(M,\Sigma_{M},\mu,T)$ is
$\varepsilon$-congruent to this Markov process.
\end{theorem}
For a proof, see subsection~\ref{PT2}. The assumptions in this
theorem are fulfilled by all Bernoulli systems in science.

The following theorem shows that, for our notion of replacement at
every observation level, also \textit{only} Bernoulli systems can
be replaced by irreducible and aperiodic Markov processes.

\begin{theorem}\label{T3}
The deterministic representation of any irreducible and aperiodic
Markov process is a Bernoulli system.
\end{theorem}
For a proof of this deep theorem, see
\citet[pp.~45--47]{Ornstein1974}.

For example, consider the baker's system $(M,\Sigma_{M},\mu,T)$
(Example~2), where $d_{M}$ is the Euclidean metric. It is a
Bernoulli system. Thus for every $\varepsilon>0$ there is a Markov process such that the baker's
system is $\varepsilon$-congruent to this Markov process. Let me
explain this. For an arbitrary $\varepsilon>0$ choose
$n\in\field{N}$ such that $\frac{\sqrt{2}}{2^{n}}<\varepsilon$.
Consider the partition
$\alpha_{n}=\{\alpha_{1},\alpha_{2},\ldots,\alpha_{2^{2n}}\}\!:=$\begin{equation}
\{[0,\frac{1}{2^{n}})\times[0,\frac{1}{2^{n}})\setminus \Gamma,
[0,\frac{1}{2^{n}})\times[\frac{1}{2^{n}}, \frac{2}{2^{n}})\setminus\Gamma,\ldots,
[\frac{2^{n}-1}{2^{n}},1]\times[\frac{2^{n}-1}{2^{n}},1]\setminus\Gamma\}.\end{equation}
Now let
$\Phi_{\alpha_{n}}(m):=\sum_{i=1}^{2^{2n}}o_{\alpha_{i}}\chi{\alpha_{i}}(m)$,
where
\begin{eqnarray}o_{\alpha_{1}}=(\frac{\sqrt{2}}{2^{n+1}},\frac{\sqrt{2}}{2^{n+1}}),
o_{\alpha_{2}}=(\frac{\sqrt{2}}{2^{n+1}},\frac{2+\sqrt{2}}{2^{n+1}}),
\ldots,\nonumber\\
o_{\alpha_{2^{2n}}}=(\frac{2^{n+1}-2+\sqrt{2}}{2^{n+1}},
\frac{2^{n+1}-2+\sqrt{2}}{2^{n+1}}).
\end{eqnarray}
It is not hard to see that
$\{\Phi_{\alpha_{n}}(T^{t});\,t\in\field{Z}\}$ is an irreducible
and aperiodic Markov process whose deterministic representation is
isomorphic to $(M,\Sigma,\mu,T)$. $(M,\Sigma,\mu,T)$ is
$\varepsilon$-congruent to this Markov process since
\begin{equation}d_{M}(m,\Phi_{\alpha_{n}}(m))\leq\frac{\sqrt{2}}{2^{n}}<\varepsilon\,\,\,\,\textnormal{for
all}\,\,m\in M.\end{equation}

Recall that irreducible and aperiodic Markov processes are widely
used in science, and they are even regarded as being second most
random. Also recall that
several deterministic systems in science are Bernoulli systems
(subsection \ref{DSRSS}). Hence Theorem~\ref{T1} and
Theorem~\ref{T2} show that irreducible and aperiodic Markov
processes are the most random stochastic processes which are
needed in order to replace deterministic systems at every
observation level. This implies that the conjecture advanced at
the beginning of this subsection is wrong: \textit{it is possible to replace measure-theoretic deterministic systems in science at every observation level by stochastic processes in science.}

\subsubsection{Previous philosophical discussion}\label{PM2}
Let me discuss the previous philosophical papers about the topic of this section. \cite{Suppes1996} and
\cite{Suppes1999} discuss an instance of Theorem~\ref{T2}, namely
that for discrete versions of billiards with convex obstacles and for any
$\varepsilon>0$ there is a Markov process such that the billiard
system is $\varepsilon$-congruent to this Markov process.
\cite{Suppes1993} (albeit with only half a page on the topic of
this section) and \cite{Winnie1998} discuss the theorem that
some continuous-time deterministic systems can be replaced at
every observation level by semi-Markov processes.

\citet[p.~196]{Suppes1996}, \citet[p.~317]{Winnie1998} and
\citet[p.~181--182]{Suppes1999} claim that the philosophical
significance of these results is that for chaotic motion and every
observation level one can choose between a deterministic
description in science and a stochastic description. For instance,
\citet[p.~196]{Suppes1996} comment on the significance of these
results:
\begin{quote}
What is fundamental is that independent of this variation of
choice of examples or experiments is that [\textit{sic}] when we do have chaotic
phenomena [...]\ then we are in a position to choose either a
deterministic or stochastic model.
\end{quote}

However, I submit that these claims are weak, and the $\varepsilon$-congruence results show more. As argued in section~\ref{DRS}, the basic results on observational
equivalence already show that for many deterministic systems, including many deterministic systems in science, the following holds: for any finite-valued observation function one can choose between a nontrivial stochastic or a deterministic description. This implies that, in a way, many deterministic systems can be replaced at every observation level by nontrivial stochastic processes. And as one would expect,
given a deterministic system in science which satisfies the assumption of Proposition~\ref{ep}, for every $\varepsilon>0$ there is a nontrivial stochastic process which is $\varepsilon$-congruent to the system. For basically all ergodic deterministic systems in science have a generating partition (Definition~\ref{generating}).\footnote{Basically all deterministic systems in science have finite Kolmogorov-Sinai entropy; and ergodic systems with finite Kolmogorov-Sinai entropy have a generating partition (cf.\ \citeauthor{Petersen1983}, \citeyear{Petersen1983}, p.~244; \citeauthor{OrnsteinWeiss1991}, \citeyear{OrnsteinWeiss1991}, p.~19).} Besides, one easily sees the following: assume that $(M,\Sigma_{M},\mu,T)$ is an ergodic measure-preserving deterministic system with a generating partition, and that $(M,d_{M})$ is separable and $\Sigma_{M}$ contains all open balls of $(M,d_{M})$. Then for every $\varepsilon>0$ there is a stochastic process $\{\Phi(T^{t});\,t\in\field{Z}\}$, where $\Phi:M\rightarrow M$ is finite-valued, which is $\varepsilon$-congruent to the
system.\footnote{Let $\varepsilon>0$. Since $(M,d_{M})$ is separable, there exists a $r\in\field{N}$ and $m_{i}\in M$, $1\leq i \leq r$, such that $\mu(M\setminus\cup_{i=1}^{r}B(m_{i},\frac{\varepsilon}{2}))<\frac{\varepsilon}{2}$ ($B(m,\varepsilon)$ is the ball of radius $\varepsilon$ around $m$). Let $\alpha$ be a generating partition. Then for each $B(m_{i},\frac{\varepsilon}{2})$ there is an $n_{i}\in\field{N}$ and a $C_{i}$ of union of elements in $\vee_{j=-n_{i}}^{n_{i}}T^{j}(\alpha)$ such that $\mu((B(m_{i},\frac{\varepsilon}{2})\setminus C_{i})\cup
(C_{i}\setminus B(m_{i},\frac{\varepsilon}{2}))<\frac{\varepsilon}{2r}$. Define $n\!:=\!\max\{n_{i}\}$, $\beta\!:=\!\{\beta_{1},\ldots,\beta_{l}\}\!:=\!\vee_{j=-n}^{n}T^{j}(\alpha)$ and  $\Phi\!:=\!\sum_{i=1}^{l}o_{i}\chi_{\beta_{i}}$ with $o_{i}\in \beta_{i}$. $\Phi$ is finite-valued and, since $\beta$ is generating, $(M,\Sigma_{M},\mu,T)$ is isomorphic via $\phi$ to the deterministic representation $(M_{2},\Sigma_{M_{2}}, \mu_{2},T_{2},\Phi_{0})$ of the stochastic process $Z_{t}:=\{\Phi(T^{t}); t\in\field{Z}\}$ \cite[p.~274]{Petersen1983}. And, by construction, $d_{M}(x,\Phi_{0}(\phi(x)))<\varepsilon$ except for a set in $M$ smaller than $\varepsilon$.}
And similar results for
chaotic systems were known long before the
$\varepsilon$-congruence results were proved
(cf.~subsection~\ref{DRS}). \textit{Hence the fact that at every
observation level one has a choice between a deterministic
description in science and a stochastic process was known long
before the $\varepsilon$-congruence results were proved, and so
cannot be the philosophical significance of these
results as claimed by these authors.}\footnote{The reader should also be warned that there
are some technical lacunae in \cite{Suppes1996} and
\cite{Suppes1999}. For instance, according to their
definition, any two systems whatsoever are $\varepsilon$-congruent
(let the metric space simply consist of one element). Also, these
authors do not seem to be aware that the continuous-time
$\varepsilon$-congruence results require the motion to be a
Bernoulli flow and so do not generally hold for ergodic systems. And
in these papers it is wrongly assumed that the notions of
isomorphism and $\varepsilon$-congruence require that the
deterministic system is looked at through an observation function (\citeauthor{Suppes1996},~\citeyear{Suppes1996}, p.~195--196,
p.~200; p.~198-200;
\citeauthor{Suppes1999},~\citeyear{Suppes1999}, p.~192, p.~195;
pp.~189--192).} As I have argued in subsection \ref{muede}, the significance of the
$\varepsilon$-congruence results is something stronger: namely
that it is possible to replace deterministic systems in science at every observation level by stochastic processes in science.

Furthermore, \cite{Suppes1996}, \cite{Winnie1998} and
\cite{Suppes1999} do not seem to be aware that \textit{also for
non-chaotic systems there is a choice between a deterministic and
a stochastic description.} To show this, I do not have to discuss
the hard question of how to define chaos. It will suffice to show
that Proposition~\ref{ep} also applies to systems which
are uncontroversially neither chaotic nor locally chaotic. Consider the measure-preserving
deterministic system $(M,\Sigma_{M},\mu,T)$ where $M:=[0,1)$
represents the unit circle, i.e.\ each $m\in\,\,M$ represents the
point $e^{2\pi m i}$, $\Sigma_{M}$ is the Lebesgue
$\sigma$-algebra, $\mu$ is the Lebesgue measure, and $T$ is the
rotation $T(m):=m+\alpha\,(\textnormal{mod}\,1)$, where
$\alpha\in\field{R}$ is irrational. It is uncontroversial that this system is neither chaotic nor locally chaotic because all solutions are stable,
i.e.\ nearby solutions stay close for all times. However, one
easily sees that it satisfies the assumption of
Proposition~\ref{ep}.\footnote{This is so because any such deterministic system with irrational $\alpha$ is ergodic (Petersen 1983, p.~49).} Consequently, this deterministic system is replaced at every observation level by a nontrivial stochastic process.\footnote{This example can be generalised: any rationally independent rotation on a torus is uncontroversially non-chaotic but fulfills the assumption of Proposition~\ref{ep} \citep[cf.][p.~51]{Petersen1983}.}

There remains the question: if one can choose between
a deterministic system and a stochastic process, which description
is preferable? \citet[pp.~317--318]{Winnie1998} dismisses \citeauthor{Suppes1993}'s (\citeyear{Suppes1993}, p.~254)
claim that in the case of the $\varepsilon$-congruence results
both descriptions are equally good.  \citeauthor{Winnie1998} argues
that the deterministic description is preferable: assume a stochastic
process replaces a deterministic system for the current
observation level. At some point in the future the observational
accuracy may be so fine that another stochastic process will be
needed to replace the deterministic system, and so on. Because there is in principle
no limitation on the observational accuracy, there is no stochastic
process that one can be sure, for practical purposes, will always give the same predictions as the deterministic system.
Hence the deterministic description is preferable.

However, I think neither \citeauthor{Winnie1998}'s (\citeyear{Winnie1998}) nor \citeauthor{Suppes1993}' (\citeyear{Suppes1993}) view is tenable. In a way, if the phenomenon under consideration is really stochastic, the stochastic description is preferable, even if the stochasticity is at a small scale and thus not observable. Likewise, if the phenomenon is really deterministic, the deterministic description is preferable. Now assume one cannot know for sure whether the phenomenon is deterministic or stochastic. Which description is then preferable in the sense of being preferable relative to our current knowledge and evidence? The answer depends on many factors, such as the kind of phenomenon under consideration, theories about fundamental physics, etc. And it may well be that the stochastic description is preferable if, for instance, a fundamental theory suggests this and one aims for a description at the most fundamental level. To sum up, neither \citeauthor{Winnie1998}'s  nor \citeauthor{Suppes1993}' view is tenable, and the question of which description is preferable needs more careful examination.

\section{Conclusion}\label{Conclusion}
The central question of this paper has been: are
deterministic and indeterministic descriptions observationally equivalent in the sense that deterministic
descriptions, when observed, and indeterministic descriptions give
the same predictions? I have tackled it for discrete-time stochastic processes and
measure-theoretic deterministic systems, both of which are
ubiquitous in science.

I have demonstrated that every stochastic process is observationally equivalent to a deterministic system, and that many deterministic systems are observationally equivalent to stochastic processes. Still, one might guess that the
measure-theoretic deterministic systems which are observationally
equivalent to stochastic processes in science do not include any
deterministic systems in science. I have shown this to be false because some deterministic systems in science even produce Bernoulli processes.
Despite this, one might guess that measure-preserving deterministic systems in science cannot give the same predictions at every observation level as stochastic processes in science. I have shown that there is indeed a limitation on observational equivalence, namely deterministic systems in science cannot give the same predictions at every observation level as Bernoulli processes. However, the guess is still wrong because I have shown (one of the $\varepsilon$-congruence results) that several deterministic systems in science give the same predictions at every observation level as Markov processes. Therefore, even kinds of stochastic processes and kinds of deterministic systems which intuitively seem to give very different predictions are observationally equivalent.

Furthermore, I have criticised the previous philosophical literature, namely \cite{Suppes1996}, \cite{Winnie1998} and \citeauthor{Suppes1999}
(\citeyear{Suppes1999}). They argue that the philosophical
significance of the $\varepsilon$-congruence results is that for
chaotic motion one can choose at every observation level between a stochastic or a deterministic description. However, this is already shown by the basic results in
subsection~\ref{DRS}. The philosophical significance of the $\varepsilon$-congruence result is really something stronger, namely, that there are deterministic systems in science that give the same predictions at every observation level as stochastic processes in science. Furthermore, these authors seem not to be aware that there are also uncontroversially non-chaotic deterministic systems which can be replaced at every observation level by stochastic processes.

\section{Appendix: Proofs}
\subsection{Proof of Proposition~\ref{ep}}\label{PP1}
\textbf{Proposition~\ref{ep}} \textit{Assume that the deterministic system $(M,\Sigma_{M},\mu,T)$ is
ergodic and that there does not exist an $n\in\field{N}$ and a
$C\in\Sigma_{M},\,\,0<\mu(C)<1,$ such that, except for
a set of measure zero, $T^{n}(C)=C$. Then for every nontrivial finite-valued
observation function $\Phi:M\rightarrow M_{O}$ and the stochastic
process $\{Z_{t}\!:=\!\Phi(T^{t});\,t\in\field{Z}\}$ the following
holds: there is an $o_{i}\in M_{O}$ such that for all
$o_{j}\in M_{O}$, $P\{Z_{t+1}\!=\!o_{j}\,|\,Z_{t}\!=\!o_{i}\}<1$.}\\\\
\textit{Proof}: I have not found a proof of this result in the
literature and thus provide one here. Notice that it suffices to
prove the following:
\begin{quote}$(*)$ Assume that $(M,\Sigma_{M},\mu,T)$ is ergodic and that it is not the case that there exists an $n\in\field{N}$ and a $C\in\Sigma_{M},$ $0<\mu(C)<1,$ such that, except for a set of measure zero, $T^{n}(C)=C$. Then for any nontrivial partition $\alpha=\{\alpha_{1},\ldots,\alpha_{n}\}$ there is an $i\in\{1,\ldots,n\}$ such that for all $j,\,1\!\leq\! j\!\leq\! n$, $\mu(T(\alpha_{i})\!\setminus\! \alpha_{j})\!>\!0$.
\end{quote}
For recall that any finite observation function has a
corresponding partition (cf.~subsection \ref{DRS}). Hence the
conclusion of $(*)$ implies that for any nontrivial finite
observation function $\Phi:M\rightarrow M_{O}$ there is an outcome
$o_{i}\in M_{O}:=\cup_{k=1}^{n}o_{k}$, $n\in\field{N}$, such that
for all possible outcomes $o_{j}\in M_{O}$ it follows that
$P\{Z_{t+1}=o_{j}\,|\,Z_{t}=o_{i}\}<1$, $t\in\field{Z}$.

So assume that the conclusion of $(*)$ does not hold, i.e.~there
exists a nontrivial partition $\alpha$ such that for each
$\alpha_{i}$ there exists an $\alpha_{j}$ with,  except for a set of measure zero,
$T(\alpha_{i})\subseteq \alpha_{j}$.

\textit{Case 1}: for all $i$ there is a $j$ such that,  except for a set of measure zero, $T(\alpha_{i})= \alpha_{j}$. Then
ergodicity implies that $\alpha_{1}$ is mapped, except for a set
of measure zero, onto all $\alpha_{k}$, $2\leq k \leq n$, before
being mapped onto itself. But this contradicts the assumption that
it is not the case that there exists an $n\in\field{N}$ and a
$C\in\Sigma_{M},$ $0<\mu(C)<1,$ such that,  except for a set of measure zero, $T^{n}(C)=C$.

\textit{Case 2}: for some $i$ there is a $j$ with,  except for a set of measure zero,
$T(\alpha_{i})\subset \alpha_{j}$ and $\mu(\alpha_{i})<\mu(\alpha_{j})$. Ergodicity implies that there exists a $k \in \field{N}$ such that,  except for a set of measure zero, $T^{k}(\alpha_{j})\subseteq \alpha_{i}$. Hence it holds that $\mu(\alpha_{j})\leq
\mu(\alpha_{i})$, yielding a contradiction, viz.\
$\mu(\alpha_{i})<\mu(\alpha_{j})\leq\mu(\alpha_{i})$.

\subsection{Proof of Theorem~\ref{T1}}\label{PT1}
\textbf{Theorem~\ref{T1}} \textit{Let $(M,\Sigma_{M},\mu,T)$ be a
deterministic system where $\Sigma_{M}$ contains all open balls of
the metric space $(M,d_{M})$, $T$ is continuous at a point $x\in
M$, every open ball around $x$ has positive measure, and there is
a set $D\in \Sigma_{M}$, $\mu(D)>0$, with
$d(T(x),D):=\inf\{d(T(x),m)\,|\,m \in D\}>0$. Then there is some
$\varepsilon>0$ for which there is no Bernoulli process to which
$(M,\Sigma_{M},\mu,T)$ is $\varepsilon$-congruent.}\\\\
\textit{Proof}: I have not found a proof of this result in the
literature and thus provide one here.

For $m\in M$, $E\subseteq M$ and $\varepsilon>0$ let the ball of
radius $\varepsilon$ around $m$ be $B(m,\varepsilon):=\{y\in
M\,|\,d(y,m)<\varepsilon\}$ and let $B(E,\varepsilon):=\cup_{m\in
E}B(m,\varepsilon)$. Since $d(T(x),D)>0$, one can choose $\gamma>0$
and $\beta>0$ such that $B(T(x),2\gamma)\cap
B(D,2\beta)=\emptyset$. Because $T$ is continuous at $x$, one can
choose $\delta>0$ such that $T(B(x,4\delta))\subseteq
B(T(x),\gamma)$. Recall that $\mu(B(x,2\delta))=\rho_{1}>0$ and
that $\mu(D)=\rho_{2}>0$. Let $\varepsilon>0$ be such that
$\varepsilon<\frac{\rho_{1}\rho_{2}}{8}$, $\varepsilon<\delta$,
$\varepsilon<\beta$ and $\varepsilon<\gamma$. I am going to show
that there is no Bernoulli process such that
$(M,\Sigma_{M},\mu,T)$ is $\varepsilon$-congruent to this
Bernoulli process.

Assume that $(M,\Sigma_{M},\mu,T)$ is $\varepsilon$-congruent to a
Bernoulli process, and let
$(\Omega,\Sigma_{\Omega},\nu,S,\Phi_{0})$ be the deterministic
representation of this Bernoulli process. This implies that
$(M,\Sigma_{M},\mu,T)$ is isomorphic (via $\phi:\hat{M}\rightarrow
\hat{\Omega}$) to the Bernoulli shift
$(\Omega,\Sigma_{\Omega},\nu,S)$ and hence that
$(M,\Sigma_{M},\mu,T)$ is a Bernoulli system. Let
$\alpha_{\Phi_{0}}:=\{\alpha_{\Phi_{0}}^{1}\ldots\alpha_{\Phi_{0}}^{s}\}$,
$s\in\field{N}$, be the partition of
$(\Omega,\Sigma_{\Omega},\nu)$ corresponding to the observation
function $\Phi_{0}$ (cf.\ subsection~\ref{DRS}). Let
$\check{M}\!:=\!M\setminus \hat{M}$ and
$\check{\Omega}\!\!:=\!\!\Omega\setminus\hat{\Omega}$. Clearly,
$\phi^{-1}(\alpha_{\Phi_{0}})\!\!:=\!\!\{\phi^{-1}(\alpha_{\Phi_{0}}^{1}\setminus
\check{\Omega})\cup\check{M},\phi^{-1}(\alpha_{\Phi_{0}}^{2}\setminus
\check{\Omega}),$  $\ldots,
\phi^{-1}(\alpha_{\Phi_{0}}^{s}\setminus \check{\Omega})\}$ is a
partition of $(M,\Sigma_{M},\mu)$.

Consider all the sets in $\phi^{-1}(\alpha_{\Phi_{0}})$ which are
assigned values in $B(x,3\delta)$, i.e.~all the sets
$a\in\phi^{-1}(\alpha_{\Phi_{0}})$ with $\Phi_{0}(\phi(m))\in
B(x,3\delta)$ for almost all $m\in a$. Denote these sets by
$A_{1},\ldots A_{n}$, $n\in\field{N}$, and let
$A:=\cup_{i=1}^{n}A_{i}$. Because $(M,\Sigma_{M},\mu,T)$ is
$\varepsilon$-congruent to
$(\Omega,\Sigma_{\Omega},\nu,S,\Phi_{0})$, it follows that
$\mu(A\setminus B(x,4\delta))<\varepsilon$ and $\mu(A\cap
B(x,2\delta))\geq \rho_{1}/2$.

Now consider all the sets in $\phi^{-1}(\alpha_{\Phi_{0}})$ which
are assigned values in $B(D,\beta)$, i.e.~all the sets
$c\in\phi^{-1}(\alpha_{\Phi_{0}})$ where $\Phi_{0}(\phi(m))\in
B(D,\beta)$ for almost all $m\in c$. Denote these sets by
$C_{1},\ldots C_{k}$, $k\in\field{N}$, and let
$C:=\cup_{i=1}^{k}C_{i}$. Because $(M,\Sigma_{M},\mu,T)$ is
$\varepsilon$-congruent to
$(\Omega,\Sigma_{\Omega},\nu,S,\Phi_{0})$, I have $\mu(C\cap
D)\geq\rho_{2}/2$ and $\mu(C\cap B(T(x),\gamma))<\varepsilon$.

Because $(\Omega,\Sigma_{\Omega},\nu,S,\Phi_{0})$ is a Bernoulli
process isomorphic to $(M,\Sigma_{M},\mu,T)$, it must hold that
$\mu(T(A_{i})\cap C_{j})=\mu(A_{i})\mu(C_{j})$ for all $i,j$,
$1\leq i \leq n$, $1\leq j\leq k$. Hence also $\mu(T(A)\cap
C)=\mu(A)\mu(C)$. But it follows that $\mu(A)\mu(C)\geq
\frac{\rho_{1}\rho_{2}}{4}$ and that $\mu(T(A)\cap C)<\varepsilon
+ \varepsilon$, and this yields the contradiction
$\frac{\rho_{1}\rho_{2}}{4}<2\varepsilon<\frac{\rho_{1}\rho_{2}}{4}$
since it was assumed that
$\varepsilon<\frac{\rho_{1}\rho_{2}}{8}$.\\

\noindent\textit{Remark 1}. I say that a stochastic process
$\{Z_{t};\,t\in\field{Z}\}$ is a \textit{generalised Bernoulli
process} if and only if (i) it takes values in an arbitrary
measurable space $(\bar{M},\Sigma_{\bar{M}})$; (ii) the random
variables  $Z_{t}$ have the same distribution for all $t$; and
(iii) $\{Z_{t};\,t \in \field{Z}\}$ are independent. A
\textit{generalised Bernoulli shift} is the deterministic
representation of a generalised Bernoulli process. A
\textit{generalised Bernoulli system} is a deterministic system
which is isomorphic to a generalised Bernoulli shift. Now it is
important to note that Theorem~1 also holds for generalised
Bernoulli processes, i.e.:\\\\ \textbf{Theorem 1*} \textit{Let
$(M,\Sigma_{M},\mu,T)$ be a deterministic system where
$\Sigma_{M}$ contains all open balls of the metric space
$(M,d_{M})$, $T$ is continuous at a point $x\in M$, every open
ball around $x$ has positive measure, and there is a set $D\in
\Sigma_{M}$, $\mu(D)>0$, with $d(T(x),D):=\inf\{d(T(x),m)\,|\,m
\in D\}>0$. Then it is not the case that for all $\varepsilon>0$
there is a generalised Bernoulli process such that
$(M,\Sigma_{M},\mu,T)$ is $\varepsilon$-congruent to this
generalised Bernoulli process.}\\\\ \textit{Proof}: the proof goes
through as above when one considers generalised Bernoulli
processes instead of Bernoulli processes, generalised Bernoulli
shifts instead of Bernoulli shifts and generalised Bernoulli
systems instead of Bernoulli systems, and one defines
$A:=\phi^{-1}(\Phi_{0}^{-1}(B(x,3\delta))\setminus
\check{\Omega})$ and
$C:=\phi^{-1}(\Phi_{0}^{-1}(B(D,\beta))\setminus \check{\Omega})$.
\linebreak[3] Clearly, because for a generalised Bernoulli process
the random variables are independent, it still holds that
$\mu(T(A)\cap C)=\mu(A)\mu(C)$ and the proof goes through as
above.

\subsection{Proof of Theorem~\ref{T2}}\label{PT2}
\textbf{Theorem~\ref{T2}} \textit{Let $(M,\Sigma_{M},\mu,T)$ be a
Bernoulli system  where the metric space $(M,d_{M})$ is separable and $\Sigma_{M}$ contains all open balls of $(M,d_{M})$. Then for any $\varepsilon>0$ there is an irreducible and aperiodic Markov process such that $(M,\Sigma_{M},\mu,T)$ is $\varepsilon$-congruent to this Markov process}.\\\\
\textit{Proof}: I have not found a proof of this result in
the literature and thus provide one here. I need the following definition.
\begin{definition}\label{generating}
A partition $\alpha$ of $(M,\Sigma_{M},\mu,T)$ is $\mathrm{generating}$ if and only if
for every $A\in\Sigma_{M}$ there is an $n\in\field{N}$ and a set $C$ of unions of elements in $\vee_{j=-n}^{n}T^{j}(\alpha)$ such that $\mu((A\setminus C)\cup(C\setminus A))<\varepsilon$ \cite[cf.][p.~244]{Petersen1983}.
\end{definition}

Per assumption, the deterministic system $(M,\Sigma_{M},\mu,T)$ is
isomorphic via $\phi:\hat{M}\rightarrow \hat{\Omega}$ to the
deterministic representation
$(\Omega,\Sigma_{\Omega},\nu,S,\Phi_{0})$ of a Bernoulli shift
with outcome space $\bar{M}$. Let
$\alpha_{\Phi_{0}}:=\{\alpha_{\Phi_{0}}^{1}\ldots\alpha_{\Phi_{0}}^{k}\}$,
$k\in\field{N}$, be the partition of
$(\Omega,\Sigma_{\Omega},\nu)$ corresponding to the observation
function $\Phi_{0}$ (cf.\ subsection~\ref{DRS}). Let
$\check{M}:=M\setminus \hat{M}$ and
$\check{\Omega}:=\Omega\setminus\hat{\Omega}$.
$\phi^{-1}(\alpha_{\Phi_{0}}):=\{\phi^{-1}(\alpha_{\Phi_{0}}^{1}\setminus
\check{\Omega})\cup\check{M},\phi^{-1}(\alpha_{\Phi_{0}}^{2}\setminus
\check{\Omega}),\ldots, \phi^{-1}(\alpha_{\Phi_{0}}^{k}\setminus
\check{\Omega})\}$ is a partition of $(M,\Sigma_{M},\mu)$. For the
partitions $\alpha=\{\alpha_{1},\ldots,\alpha_{n}\}$ and
$\beta=\{\beta_{1},\ldots,\beta_{n}\}, n\in\field{N}$, $\alpha
\vee \beta$ is the partition $\bigcup_{1\leq i,j\leq
n}\alpha_{i}\cap\beta_{j}$. Clearly, if $\alpha$ is a partition of
$M$, $T^{-t}\alpha:=\{T^{-t}\alpha_{1},\ldots,T^{-t}\alpha_{n}\}$,
$t\in\field{Z}$, are partitions.

Since $(M,d_{M})$ is separable, there exists an $r\in\field{N}$ and
$m_{i}\in M$, \linebreak $1\leq i \leq r$, such that
$\mu(M\setminus
\cup_{i=1}^{r}B(m_{i},\frac{\varepsilon}{2}))<\frac{\varepsilon}{2}$.
Because for a Bernoulli system $\Phi^{-1}(\alpha_{\Phi_{0}})$ is generating \citep[p.~275]{Petersen1983}, for each $B(m_{i},
\frac{\varepsilon}{2})$ there is an $n_{i}\in\field{N}$ and a
$C_{i}$ of union of elements in
$\vee_{j=-n_{i}}^{n_{i}}T^{j}(\phi^{-1}(\alpha_{\Phi_{0}}))$ such
that $\mu(D_{i})<\frac{\varepsilon}{2r}$, where
$D_{i}:=(B(m_{i},\frac{\varepsilon}{2})\setminus C_{i})\cup
(C_{i}\setminus B(m_{i},\frac{\varepsilon}{2}))$. Define
$n:=\max\{n_{i}\}$. For
$Q=\{q_{1},\ldots,q_{l}\}:=\vee_{j=-n}^{n}S^{j}(\alpha_{\Phi_{0}})$
let $\Phi_{0}^{Q}:\Omega\rightarrow M;
\Phi_{0}^{Q}(\omega)=\sum_{i=1}^{l}o_{i}\chi_{q_{i}}(\omega)$,
where $o_{i}\in \phi^{-1}(q_{i}\setminus \check{\Omega})$. Note
that $o_{i}\neq o_{j}$ for $i\neq j$, $1\leq i,j\leq l$. Then
\begin{equation}\label{alpha}
d_{M}(m,\Phi^{Q}_{0}(\phi(m)))<\varepsilon\,\,\textnormal{except
for a set in $M$ of measure}\,\,<\varepsilon.
\end{equation}

Now let $(X,\Sigma_{X},\lambda,R,\Theta_{0})$ be the deterministic
representation of the stochastic process
$\{\Phi_{0}^{Q}(S^{t});\,t\in\field{Z}\}$ from
$(\Omega,\Sigma_{\Omega},\nu)$ to $(M,\Sigma_{M})$. This process
is a Markov process since for any $k\in\field{N}$ and any
$A,B_{1},\ldots,B_{k} \in \bar{M}^{2n+1}$,
\begin{eqnarray}
&&
  \frac{\nu(\{\omega\in\Omega\,|\,(\omega_{-n}\ldots\omega_{n}) = A\,\,\,\textnormal{and}\,\,\,(\omega_{-n+1}\ldots\omega_{n+1}) = B_{1}\})}{
 \nu(\{\omega\in\Omega\,|\,(\omega_{-n+1}\ldots\omega_{n+1}) =B_{1}\})
}=\\ &&
\frac{\nu(\{\omega\!\in\!\Omega|(\omega_{-n}\!\ldots\!\omega_{n})\!\!=\!\!
A\,\textnormal{and}(\omega_{-n+1}\!\ldots\!\omega_{n+1})\!\!=\!\!
B_{1},  \!\ldots\!,(\omega_{-n+k}\!\ldots\!\omega_{n+k})\!\!=\!\!
B_{k}\})}{\nu(\{\omega\!\in\!\Omega|(\omega_{-n+1}\!\ldots\!\omega_{n+1})\!\!=\!\!
B_{1},\!\ldots\!,(\omega_{-n+k}\!\ldots\!\omega_{n+k})\!\!=\!\!
B_{k}\})\nonumber},
\end{eqnarray}
\begin{eqnarray}
\textnormal{if}\,\,\nu(\{\omega\!\in\!\Omega|(\omega_{-n}\!\ldots\!\omega_{n})\!\!=\!\!
A\,\textnormal{and}(\omega_{-n+1}\!\ldots\!\omega_{n+1})\!\!=\!\!
B_{1},  \!\ldots\!,(\omega_{-n+k}\!\ldots\!\omega_{n+k})\!\!=\!\!
B_{k}\})\!\!>\!\!0.\nonumber
\end{eqnarray}

Because $S$ is a shift, one sees that for all $i,j$, $1\leq
i,j\leq l$, there is a $k\geq 1$ such that $P^{k}(o_{i},o_{j})>0$,
and hence that the Markov process is irreducible. One also sees
that there exists an outcome $o_{i}$, $1\leq i\leq l$, such that
$P^{1}(o_{i},o_{i})>0$. Hence $d_{o_{i}}=1$; and since all
outcomes of an irreducible Markov process have the same
periodicity \cite[p.~131]{Cinlar1975}, it follows that the Markov
process is also aperiodic.

Consider $\psi:\Omega\rightarrow X$,
$\psi(\omega)=\ldots\Phi^{Q}_{0}(S^{-1}(\omega)),\Phi^{Q}_{0}(\omega),\Phi^{Q}_{0}(S(\omega))\ldots$,
\linebreak for $\omega \in \Omega$. Clearly, there is a set
$\hat{X}\subseteq X$ with $\lambda(\hat{X})=1$ such that
$\psi:\Omega\rightarrow\hat{X}$ is bijective and
measure-preserving and $R(\psi(\omega))=\psi(S(\omega))$ for all
$\omega\in\Omega$. Hence $(\Omega,\Sigma_{\Omega},\nu,S)$ is
isomorphic to $(X,\Sigma_{X},\lambda,R)$ via $\psi$, and thus
$(M,\Sigma_{M},\mu,T)$ is isomorphic to $(X,\Sigma_{X},\lambda,R)$
via $\theta=\psi(\phi)$. Now because of (\ref{alpha}) it holds
that
\begin{equation}
d_{M}(m,\Theta_{0}(\theta(m)))<\varepsilon\,\,\textnormal{except
for a set in $M$ of measure}\,<\varepsilon.
\end{equation}

\section*{Acknowledgments} I am indebted to Jeremy Butterfield for his continued support and, in particular, for getting up at 3 a.m.\ to read my papers. For valuable comments I also want to thank Cymra Haskell, Franz Huber, Thomas M\"{u}ller, Donald Ornstein, Amy Radunskaya, Peter Smith, Jos Uffink, two anonymous referees, and the audiences at the Oxford Philosophy of Physics Discussion Group, the Cambridge Moral Sciences Club, and the Wittgenstein Symposium 2008. I am grateful to St John's College Cambridge for financial support.

\end{document}